\documentclass[12pt,a4paper,oneside]{amsart}

\usepackage[utf8]{inputenc}
\usepackage{amssymb}
\usepackage{amsthm}
\usepackage{bm} 
\usepackage[top=1.7cm, bottom=1.7cm, left=2.7cm, right=1.7cm]{geometry}
\usepackage[shortlabels]{enumitem}   
\setlist[itemize]{noitemsep,topsep=2pt,leftmargin=16pt}
\setlist[enumerate]{noitemsep,topsep=2pt,leftmargin=16pt} 

\usepackage{makecell}  

\newtheorem{Thm}{Theorem}

\newtheorem{Rem}{Remark}
\newenvironment{Pf}{ Proof.}{\(\square\)}

\title[Randers metrics with compatible linear connections...]{Randers metrics with compatible linear connections: a coordinate-free approach}
\author{Márk Oláh}
\address{Institute of Mathematics, University of Debrecen, H-4002 Debrecen, P.O.Box 400, Hungary}
\email{olah.mark@science.unideb.hu}
\author{Csaba Vincze}
\address{Institute of Mathematics, University of Debrecen, H-4002 Debrecen, P.O.Box 400, Hungary}
\email{csvincze@science.unideb.hu}

\begin{document}

\begin{abstract} A Randers space is a differentiable manifold equipped with a Randers metric. It is the sum of a Riemannian metric and a one-form on the base manifold. The compatibility of a linear connection with the metric means that the parallel transports preserve the Randers norm of tangent vectors. The existence of such a linear connection is not guaranteed in general. If it does exist then we speak about a generalized Berwald Randers metric. In what follows we give a necessary and sufficient condition for a Randers metric to be a generalized Berwald metric and we describe some distinguished compatible linear connections. The method is based on the solution of constrained optimization problems for tensors that are in one-to-one correspondence to the compatible linear connections. The solutions are given in terms of explicit formulas by choosing the free tensor components to be zero. Throughout the paper we use a coordinate-free approach to keep the geometric feature of the argumentation as far as possible.  
\end{abstract}

\maketitle

\section{Introduction}

Let $M$ be a connected differentiable manifold equipped with a Riemannian metric $\alpha$ and consider a one-form $\beta$ on the base manifold such that the Riemannian norm of its dual vector field $\beta^{\sharp}$ is less than one. The sum of the Riemannian norm and such a one-form gives a \emph{Randers metric}  
$$F\colon TM\to \mathbb{R}, \quad F(v):=\sqrt{\alpha(v,v)}+\beta(v)$$ 
on the base manifold. Since the Riemannian norm of the dual vector field is less than one, the origin is in the  interior of the unit balls. A linear connection $\nabla$ on $M$ is \emph{compatible} with $F$ if the parallel transports with respect to $\nabla$ preserve the Randers norm of tangent vectors. Using that $\beta$ is a linear functional at a given point, a linear parallel translation gives that
$$\sqrt{\alpha(\tau(v), \tau(v))}+\beta\circ \tau(v)=\sqrt{\alpha(v,v)}+\beta(v)$$
and
$$\sqrt{\alpha(\tau(v), \tau(v))}-\beta\circ \tau(v)=\sqrt{\alpha(v,v)}-\beta(v)$$
by changing the role of $v$ and $-v$. Adding/Subtracting the formulas, it follows that
$$\sqrt{\alpha(\tau(v), \tau(v))}=\sqrt{\alpha(v,v)}\quad \textrm{and} \quad \beta\circ \tau(v)=\beta(v)$$
and we have the following result.
\begin{Thm} {\emph{\cite{Vin1}}} A linear connection $\nabla$ is compatible with the Randers metric if and only if 
$\nabla \alpha=0$ and $\nabla \beta=0$.
\end{Thm}

Taking the musical isomorphisms as Riemannian isometries we can introduce an inner product, a norm etc. for tensors at a given point of the base manifold \cite{Lee} (Proposition 2. 40):  if $e_1$, $\ldots$, $e_n\in T_pM$ is an orthonormal basis and $\theta^1$, $\ldots$, $\theta^n$ is the corresponding dual basis, then the set
$$e_{i_1}\otimes \ldots \otimes e_{i_k}\otimes \theta^{j_1}\otimes \ldots \otimes \theta^{j_l} \quad (i_1, \ldots, i_k, j_1, \ldots, j_l=1, \ldots, n)$$
forms an orthonormal basis for the vector space of tensor fields at the point $p\in M$. Having an orthonormal basis, the inner product can be computed as the usual dot product in terms of the tensor components. 

\section{An optimization problem for the difference tensor} 
Let $\nabla^*$ be the L\'{e}vi-Civita connection of $\alpha$. Following the L\'{e}vi-Civita process, any metric linear connections satisfying $\nabla \alpha=0$ can be given in terms of the difference tensor
\begin{equation}
\nabla_X Y=\nabla^*_X Y+A(X,Y),
\end{equation}
where
\begin{equation}
\label{cond:01}
\alpha(A(X,Y),Z)=-\alpha(A(X,Z), Y)
\end{equation}
because of $\nabla \alpha=0$. In addition, condition $\nabla \beta=0$ is equivalent to
\begin{equation}
\label{cond:02}
A(X,\beta^{\sharp})=-\nabla^*_X \beta^{\sharp}.
\end{equation}
Introducing the inclusion operator $\iota_X A \ (Y):=A(X,Y)$, we are looking for the pointwise solution of the constrained optimization problem
\begin{equation}
\label{optprob:01}
\textrm{min\ } \left | \iota_X A\right |^2 \quad \textrm{subject to \ } \alpha(A(X,Y),Z)=-\alpha(A(X,Z), Y) \quad \textrm{and} \quad A(X,\beta^{\sharp})=-\nabla^*_X \beta^{\sharp}
\end{equation}
for any vector field $X$.

\begin{Thm}
\label{solvability:01}
The constrained optimization problem \eqref{optprob:01} is solvable for any vector field $X$ if and only if $\beta^{\sharp}$ is of constant Riemannian length.
\end{Thm}

\begin{Pf}
If the problem \eqref{optprob:01} is solvable for any vector field $X$, then we have that
\begin{equation}
\label{cond:03}
0=\alpha(A(X,\beta^{\sharp}), \beta^{\sharp})=\alpha(\nabla^*_X \beta^{\sharp}, \beta^{\sharp})
\end{equation}
that is $\beta^{\sharp}$ is of constant length $\alpha(\beta^{\sharp}, \beta^{\sharp})=K^2$. Conversely, if $K=0$ then the optimization problem is solved by the identically zero difference tensor. Otherwise, if $K\neq  0$ is a positive constant then we can write the vector fields into the form
$$Y=Y^{\bot}+Y^{\beta}, \quad \textrm{where} \quad Y^{\beta}=\frac{\alpha(Y, \beta^{\sharp})}{\ K^2}\beta^{\sharp}.$$
To minimize its contribution to the norm of the mapping $\iota_X A$, the projected part is chosen to be zero, that is 
\begin{equation}
\label{difftensor:def:01}
 \alpha(A(X, Y^{\bot}), Z^{\bot})=0 \quad \textrm{and} \quad A(X, Y^{\beta})\stackrel{\eqref{cond:02}}{=}-\frac{\alpha(Y, \beta^{\sharp})}{\ K^2}\nabla^*_X \beta^{\sharp}.
\end{equation}
Then we have that
$$ A(X,Y)=A(X, Y^{\bot})+A(X, Y^{\beta})\stackrel{\eqref{difftensor:def:01}}{=}\frac{\alpha(A(X, Y^{\bot}), \beta^{\sharp})}{K^2}\beta^{\sharp}-\frac{\alpha(Y, \beta^{\sharp})}{\ K^2}\nabla^*_X \beta^{\sharp}=$$
$$-\frac{\alpha(A(X, \beta^{\sharp}), Y^{\bot})}{K^2}\beta^{\sharp}-\frac{\alpha(Y, \beta^{\sharp})}{\ K^2}\nabla^*_X \beta^{\sharp}=
\frac{\alpha(\nabla^*_X \beta^{\sharp}, Y^{\bot})}{K^2}\beta^{\sharp}-\frac{\alpha(Y, \beta^{\sharp})}{\ K^2}\nabla^*_X \beta^{\sharp}=$$
$$\frac{\alpha(\nabla^*_X \beta^{\sharp}, Y)}{K^2}\beta^{\sharp}-\frac{\alpha(Y, \beta^{\sharp})}{\ K^2}\nabla^*_X \beta^{\sharp},$$
where 
$$\alpha(\nabla^*_X \beta^{\sharp}, Y^{\bot})=\alpha(\nabla^*_X \beta^{\sharp}, Y^{\bot}+Y^{\beta})=\alpha(\nabla^*_X \beta^{\sharp}, Y)$$
because of the constant length of the dual vector field. The vanishing of the projected component shows that the difference tensor given by \eqref{difftensor:def:01} is the (pointwise) solution of the optimization problem \eqref{optprob:01} for any vector field $X$.
\end{Pf}

If $K=0$ then $F=\alpha$ is a Riemannian metric and the compatible linear connections are the metric linear connections satisfying $\nabla \alpha=0$. If $K\neq 0$ is a positive constant then Theorem \ref{solvability:01} reproduces the criteria in \cite{Vin1} for a  Randers metric to be a generalized Berwald metric together with clarifying the geometric meaning of the construction of a special compatible linear connection. The results are summarized in the following theorem.

\begin{Thm} A non-Riemannian Randers metric admits a compatible linear connection if and only if $\beta^{\sharp}$
 is of constant Riemannian length $K\neq 0$. Then the pointwise solution of the optimization problem \eqref{optprob:01} gives a difference tensor 
\begin{equation}
\label{difftensor:def:02}
A(X,Y)=\frac{\alpha(\nabla^*_X \beta^{\sharp}, Y)}{K^2}\beta^{\sharp}-\frac{\alpha(Y, \beta^{\sharp})}{\ K^2}\nabla^*_X \beta^{\sharp}
\end{equation}
such that $\nabla^{\circ}_X Y=\nabla^*_X Y+A(X,Y)$ is a compatible linear connection with the Randers metric.
\end{Thm}

\begin{Rem} {\emph{As a topological obstruction, we can conclude that a connected differentiable manifold (satisfying the second axiom of countability) admits a non-Riemannian generalized Berwald Randers metric if and only if there exists a nowhere vanishing vector field playing the role of the dual vector field of the perturbating term with respect to a Riemannian metric.}} 
\end{Rem}
     
\section{An optimization problem for the torsion tensor} The idea of the so-called extremal compatible linear connection minimizing the norm of the torsion tensor can be found in \cite{V14}. It is a direct generalization of the L\'{e}vi-Civita connection but the torsion tensor is not necessarily zero. In case of Randers metrics, the local description is given in \cite{VinOl}. In what follows we give a coordinate-free description of the extremal compatible linear connection in case of Randers metrics provided that the dual vector field $\beta^{\sharp}$ is of constant Riemannian length $K\neq 0$. First of all let us consider the compatible linear connections in the general form 
$$\nabla_XY=\nabla^*_X Y+A(X,Y)+B(X,Y)=\nabla^{\circ}_X Y+B(X,Y),$$
where $A(X,Y)$ is defined by \eqref{difftensor:def:02}, 
\begin{equation}
\label{cond:B}
\alpha(B(X,Y),Z)=-\alpha(B(X,Z),Y) \quad \textrm{and} \quad B(X, \beta^{\sharp})=0.
\end{equation}
The torsion tensor is
$$T(X,Y)=A(X,Y)-A(Y,X)+B(X,Y)-B(Y,X)=T^{\circ}(X,Y)+B(X,Y)-B(Y,X),$$
where
\begin{equation}
\label{tor:01}
 T^{\circ}(X,Y)=\frac{\alpha(\nabla^*_X \beta^{\sharp}, Y)}{K^2}\beta^{\sharp}-\frac{\alpha(Y, \beta^{\sharp})}{\ K^2}\nabla^*_X \beta^{\sharp}-\frac{\alpha(\nabla^*_Y \beta^{\sharp}, X)}{K^2}\beta^{\sharp}+\frac{\alpha(X, \beta^{\sharp})}{\ K^2}\nabla^*_Y \beta^{\sharp}.
\end{equation}
We are going to solve the optimization problem 
\begin{equation}
\label{optprob:02}
\textrm{min\ } \left | T \right|^2 \quad \textrm{subject to \ } \alpha(B(X,Y),Z)=-\alpha(B(X,Z),Y) \quad \textrm{and} \quad B(X,\beta^{\sharp})=0.
\end{equation}
To minimize its contribution to the norm of the torsion tensor, the projected part is chosen to be zero, that is 
\begin{equation}
\label{cond:T01}
\alpha(B(X^{\bot}, Y^{\bot}), Z^{\bot})=0 \ \Rightarrow\ \alpha(T(X^{\bot}, Y^{\bot}), Z^{\bot})=\alpha(T^{\circ}(X^{\bot}, Y^{\bot}), Z^{\bot})\stackrel{\eqref{tor:01}}{=}0.
\end{equation}
Since
$$\alpha(B(X^{\bot}, Y^{\bot}), \beta^{\sharp})=\alpha(B(Y^{\bot}, X^{\bot}), \beta^{\sharp})\stackrel{\eqref{cond:B}}{=}0,$$
there are no optional components and we have that
\begin{gather}
\label{cond:T02}
\alpha(T(X^{\bot}, Y^{\bot}), \beta^{\sharp})=\alpha(T^{\circ}(X^{\bot}, Y^{\bot}), \beta^{\sharp})\stackrel{\eqref{tor:01}}{=} \alpha(\nabla^*_{X^{\bot}} \beta^{\sharp}, Y^{\bot})-\alpha(\nabla^*_{Y^{\bot}} \beta^{\sharp}, X^{\bot}).
\end{gather}
Therefore
$$T(X^{\bot}, Y^{\bot})=T^{\circ}(X^{\bot}, Y^{\bot}).$$
Two possible cases still need to investigate. Since
$$\alpha(B(X^{\bot}, \beta^{\sharp}), \beta^{\sharp})=\alpha(B(\beta^{\sharp}, X^{\bot}), \beta^{\sharp})\stackrel{\eqref{cond:B}}{=}0,$$
there are no optional components and we have that  
\begin{gather}
\label{cond:T03}
\alpha(T(X^{\bot}, \beta^{\sharp}), \beta^{\sharp})=\alpha(T^{\circ}(X^{\bot}, \beta^{\sharp}), \beta^{\sharp})\stackrel{\eqref{tor:01}}{=} -\alpha(\nabla^*_{\beta^{\sharp}} \beta^{\sharp}, X^{\bot})=\alpha(\nabla^*_{\beta^{\sharp}} X^{\bot}, \beta^{\sharp})
\end{gather}
because the dual vector field is of constant Riemannian length and $\nabla^* \alpha=0$. Finally, 
$$
\alpha(T(X^{\bot}, \beta^{\sharp}), Y^{\bot})=\alpha(T^{\circ}(X^{\bot}, \beta^{\sharp}), Y^{\bot})+\alpha(B(X^{\bot}, \beta^{\sharp}), Y^{\bot})-\alpha(B(\beta^{\sharp},X^{\bot}), Y^{\bot})\stackrel{\eqref{cond:B}}{=}$$
$$\alpha(T^{\circ}(X^{\bot}, \beta^{\sharp}), Y^{\bot})-\alpha(B(\beta^{\sharp},X^{\bot}), Y^{\bot}).$$
Changing the role of $X$ and $Y$, we have that 
$$\alpha(T(X^{\bot}, \beta^{\sharp}), Y^{\bot})+\alpha(T(Y^{\bot}, \beta^{\sharp}), X^{\bot})\stackrel{\eqref{cond:B}}{=}\alpha(T^{\circ}(X^{\bot}, \beta^{\sharp}), Y^{\bot})+\alpha(T^{\circ}(Y^{\bot}, \beta^{\sharp}), X^{\bot}).$$
Therefore, solving the constrained optimization problem
$$\min \ x^2+y^2 \ \ \textrm{subject to} \ \ x+y=\textrm{const.},$$
the contribution of the sum of type
$$\alpha^2(T(X^{\bot}, \beta^{\sharp}), Y^{\bot})+\alpha^2(T(Y^{\bot}, \beta^{\sharp}), X^{\bot})$$
to the norm of the torsion tensor is minimal if and only if 
\begin{gather}
\label{cond:T04}
\alpha(T(X^{\bot}, \beta^{\sharp}), Y^{\bot})=\alpha(T(Y^{\bot}, \beta^{\sharp}), X^{\bot})=\\ \notag
\frac{\alpha(T^{\circ}(X^{\bot}, \beta^{\sharp}), Y^{\bot})+\alpha(T^{\circ}(Y^{\bot}, \beta^{\sharp}), X^{\bot})}{2}\stackrel{\eqref{tor:01}}{=}-\frac{\alpha(\nabla^*_{X^{\bot}} \beta^{\sharp}, Y^{\bot})+\alpha(X^{\bot}, \nabla^*_{Y^{\bot}} \beta^{\sharp})}{2}.
\end{gather}

\begin{Thm}
The extremal compatible linear connection of a generalized Berwald Randers metric is
$\nabla^{\circ}$ if and only if the projected vector fields form an integrable distribution. 
\end{Thm}

\begin{Pf} The extremal compatible linear connection of a generalized Berwald Randers metric is
$\nabla^{\circ}$ if and only if $T=T^{\circ}$. Using formula \eqref{cond:T04}, 
$$\alpha(T^{\circ}(X^{\bot}, \beta^{\sharp}), Y^{\bot})=\alpha(T^{\circ}(Y^{\bot}, \beta^{\sharp}), X^{\bot}) \quad \Leftrightarrow \quad \alpha(\nabla^*_{X^{\bot}} \beta^{\sharp}, Y^{\bot})=\alpha(X^{\bot}, \nabla^*_{Y^{\bot}} \beta^{\sharp}).$$
Since $\nabla^*$ is the L\'evi-Civita connection of the Riemannian metric, it is a torsion-free metric linear connection. Therefore 
$$0=\alpha(\nabla^*_{X^{\bot}} \beta^{\sharp}, Y^{\bot})-\alpha(X^{\bot}, \nabla^*_{Y^{\bot}} \beta^{\sharp})=-\alpha (\beta^{\sharp}, \nabla^*_{X^{\bot}} Y^{\bot}-\nabla^*_{Y^{\bot}} X^{\bot})=$$
$$-\alpha (\beta^{\sharp}, [X^{\bot}, Y^{\bot}]) \ \Leftrightarrow \  [X^{\bot}, Y^{\bot}]^{\bot}= [X^{\bot}, Y^{\bot}]$$
as was to be proved.
\end{Pf}

\begin{Thm}
The torsion of the extremal compatible linear connection of a generalized Berwald Randers metric is
$$\alpha(T(X,Y),Z)=\alpha(T^{\circ}(X,Y),Z)+\frac{\beta(Y)d\beta(X,Z)-\beta(X)d\beta(Y,Z)}{2K^2}+$$
$$\frac{\beta(X)d\beta(Y, \beta^{\sharp})-\beta(Y)d\beta(X, \beta^{\sharp})}{2K^4}\beta(Z).$$
\end{Thm}

\begin{Pf}  Evaluating the term
$$\alpha(\Omega(X, Y),Z):=\frac{\beta(Y)d\beta(X,Z) -\beta(X)d\beta(Y,Z)}{2K^2}+\frac{\beta(X)d\beta(Y, \beta^{\sharp})-\beta(Y)d\beta(X, \beta^{\sharp})}{2K^4}\beta(Z)$$
at the triplets belonging to \eqref{cond:T01} - \eqref{cond:T04} we have that
$$\alpha(\Omega (X^{\bot}, Y^{\bot}), Z^{\bot})=0, \quad \alpha(\Omega (X^{\bot}, Y^{\bot}), \beta^{\sharp})=0 \quad \Rightarrow \quad \Omega (X^{\bot}, Y^{\bot})=0,$$
that is
$$T^{\circ}(X^{\bot}, Y^{\bot})+\Omega (X^{\bot}, Y^{\bot})=T^{\circ}(X^{\bot}, Y^{\bot})\stackrel{\eqref{cond:T01}, \eqref{cond:T02}}{=}T(X^{\bot}, Y^{\bot}).$$
On the other hand
$$\alpha(\Omega (X^{\bot}, \beta^{\sharp}), \beta^{\sharp})=0,$$
that is
$$\alpha(T^{\circ}(X^{\bot}, \beta^{\sharp}), \beta^{\sharp})+\alpha(\Omega (X^{\bot}, \beta^{\sharp}), \beta^{\sharp})=\alpha(T^{\circ}(X^{\bot}, \beta^{\sharp}), \beta^{\sharp})\stackrel{\eqref{cond:T03}}{=}\alpha(T(X^{\bot}, \beta^{\sharp}), \beta^{\sharp}).$$
Finally,
$$\alpha(T^{\circ}(X^{\bot}, \beta^{\sharp}), Y^{\bot})+\alpha(\Omega (X^{\bot}, \beta^{\sharp}), Y^{\bot})\stackrel{\eqref{tor:01}}{=}-\alpha(\nabla^*_{X^{\bot}} \beta^{\sharp}, Y^{\bot})+\frac{1}{2}d\beta(X^{\bot}, Y^{\bot})=$$
$$-\alpha(\nabla^*_{X^{\bot}} \beta^{\sharp}, Y^{\bot})+\frac{\alpha(\nabla^*_{X^{\bot}} \beta^{\sharp}, Y^{\bot})-\alpha(X^{\bot}, \nabla^*_{Y^{\bot}} \beta^{\sharp})}{2}=$$
$$-\frac{\alpha(\nabla^*_{X^{\bot}} \beta^{\sharp}, Y^{\bot})+\alpha(X^{\bot}, \nabla^*_{Y^{\bot}} \beta^{\sharp})}{2}\stackrel{\eqref{cond:T04}}{=}\alpha(T(X^{\bot}, \beta^{\sharp}), Y^{\bot})$$
as  was to be proved.
\end{Pf}

\begin{Rem}{\emph{Using an orthonormal system of coordinate vector fields
$$\partial_1, \ldots, \partial_{n-1}, \partial_n=\beta^{\sharp}/K$$
at the point $p\in M$, the local components of the torsion tensor of the extremal compatible linear connection can be computed by formulas \eqref{cond:T01} - \eqref{cond:T04}: for any $a, b, c=1, \ldots, n-1$
$$\alpha(T(\partial_a, \partial_b), \partial_c)\stackrel{\eqref{cond:T01}}{=}0 \ \Rightarrow \ T_{ab}^c=0,$$
$$\alpha(T(\partial_a, \partial_b), \beta^{\sharp})\stackrel{\eqref{cond:T02}}{=}\alpha(\nabla^*_{\partial_a} \beta^{\sharp}, \partial_b)-\alpha(\nabla^*_{\partial_b} \beta^{\sharp}, \partial_a)=$$
$$ \partial_a \beta_b -\partial_b \beta_a -\alpha(\beta^{\sharp}, \nabla^*_{\partial_a} \partial_b-\nabla^*_{\partial_b} \partial_a)=\partial_a \beta_b -\partial_b \beta_a \ \Rightarrow \ T_{ab}^n=\frac{1}{K}\left(\partial_a \beta_b- \partial_b \beta_a\right)$$
because the torsion of $\nabla^*$ is zero and $\partial_n=\beta^{\sharp}/K$ at the point $p\in M$. Using formula \eqref{cond:T03}, 
$$\alpha(T(\partial_a, \beta^{\sharp}), \beta^{\sharp}) \stackrel{\eqref{cond:T03}}{=} \alpha(\nabla^*_{\beta^{\sharp}} \partial_a^{\bot}, \beta^{\sharp}),$$
where
$$\partial_a^{\bot}=\partial_a -\frac{\beta_a}{K^2} \beta^{\sharp} \ \Rightarrow \ \nabla^*_{\beta^{\sharp}} \partial_a^{\bot}=\nabla^*_{\beta^{\sharp}} \partial_a -\partial_n \beta_a \partial_n-\frac{\beta_a}{K^2}\nabla^*_{\beta^{\sharp}} \beta^{\sharp}=\nabla^*_{\beta^{\sharp}} \partial_a -\partial_n \beta_a \partial_n$$
because $\partial_n=\beta^{\sharp}/K$ and $\beta_1=\ldots =\beta_{n-1}=0$ at the point $p\in M$. Therefore
$$T_{a n}^n={\Gamma^*}_{an}^n-\frac{1}{K} \partial_n \beta_a.$$
Finally,
$$\alpha(T(\partial_a, \beta^{\sharp}), \partial_c)\stackrel{\eqref{cond:T04}}{=}-\frac{\alpha(\nabla^*_{\partial_a} \beta^{\sharp}, \partial_c)+\alpha(\partial_a, \nabla^*_{\partial_c} \beta^{\sharp})}{2},$$
where
$$\alpha(\nabla^*_{\partial_a} \beta^{\sharp}, \partial_c)+\alpha(\partial_a, \nabla^*_{\partial_c} \beta^{\sharp})=\partial_a \beta_c +\partial_c \beta_a-\alpha(\beta^{\sharp}, \nabla^*_{\partial_a} \partial_c+\nabla^*_{\partial_c} \partial_a)=$$
$$\partial_a \beta_c +\partial_c \beta_a-2\alpha(\beta^{\sharp}, \nabla^*_{\partial_a} \partial_c)$$
because the torsion of $\nabla^*$ is zero. Therefore
$$T_{an}^c={\Gamma^*}_{ac}^n-\frac{1}{2K} \left(\partial_a \beta_c+\partial_c \beta_a\right)$$
because $\partial_n=\beta^{\sharp}/K$ at the point $p\in M$.}}
\end{Rem}

\section{Acknowledgement}

M\'ark Ol\'ah has received funding from the HUN-REN Hungarian Research Network.

\end{document}